\documentclass[11pt]{amsart}
\usepackage[T1]{fontenc}
\usepackage[american]{babel}
\usepackage{csquotes,enumerate,amssymb}
\usepackage[backend=bibtex,firstinits=true, url=false, isbn=false]{biblatex}
\numberwithin{equation}{section}
\newcommand{\cmd}[1]{\texttt{\textbackslash #1}}
\title{ O\lowercase{n}  p\lowercase{robabilistic} p\lowercase{roofs of} c\lowercase{ertain} b\lowercase{inomial} i\lowercase{dentities}}

\author{P. Vellaisamy}
\address{P. Vellaisamy, Department of Mathematics,
 Indian Institute of Technology Bombay, Powai, Mumbai 400076, INDIA.}
 \email{pv@math.iitb.ac.in}
\date{\today}

\begin{document}
\begin{abstract}
We give a simple statistical proof of a binomial identity, by evaluating the Laplace transform of the maximum of $n$ independent exponential random variables in two different ways. As a by product, we obtain a simple proof of an interesting result concerning the exponential distribution. The connections between a probabilistic approach and the statistical approach are discussed, which explains why certain binomial identities admit probabilistic interpretations. In the process, several new binomial identities are also obtained
and discussed.
\end{abstract}
\maketitle
\section{Introduction} 
The objectives of this article are manifold:(a) to provide a simple statistical proof of a basic binomial identity, (b) to provide a rigorous proof of an interesting result concerning the maximum of $n$ independent exponential random variables, (c) to discuss two approaches  for obtaining certain binomial identities and (d) to explain the rationale on why certain binomial identities admit probabilistic interpretations/proofs, while others may not.

Recently, Peterson \cite{peter} gave an interesting probabilistic  proof of the binomial identity
\begin{equation}\label{identity1}
 \sum_{k=0}^{n}(-1)^{k}\binom{n}{k}\left(\frac{s}{s+k}\right)=\prod_{k=1}^{n}\left(\frac{k}{s+k}\right),
\end{equation}
for $s>0$, and $n\in\mathbb{N}=\{1,2,\ldots\}$,
which we call the basic binomial identity.
He mentioned that mathematical proofs require a rather advanced tools of mathematics such as hypergeometric functions, Chu-Vandermonde formula \cite{gasper} for the right side of \eqref{identity1} or the use of Rice-integral formulas \cite{flaj} for the left side of \eqref{identity1}.
His proof involves the computation of the probability of a certain event, concerning $n$ independent exponential random variables, in two different ways and equating them. We call this the probabilistic approach. However, the probability computation in one way requires an interesting result \cite[Lemma 1]{peter} concerning the maximum of $n$ independent exponential random variables. He provided a rather heuristic argument of this result.
While going through this probabilistic proof, a natural question arises if one could provide a proof of the binomial identity in (\ref{identity1}), which is simple and does not use Lemma 1 of \cite{peter}. The main of this note is to provide such a direct statistical proof, based on Laplace transforms, and bring out its relation to the probabilistic approach used in \cite{peter}. An interesting by-product of this approach is that it readily provides a simple, yet rigorous, proof of Lemma 1 of \cite{peter}. Some new binomial identities are also obtained. Finally, we discuss an interesting relationship among binomial identities, probabilistic interpretations and the associated Laplace transforms.

\vspace{0.1cm}

\section{A Probabilistic Approach}
In this section, we briefly describe the probabilistic approach due to Peterson \cite{peter}. Let $X_{1}.\ldots,X_{n}$ be $n$ independent exponential $Exp(1)$ random variables with density 
\begin{equation*}
 f(x)=e^{-x},~~x>0,
\end{equation*}
and $X_{(n)}=\max(X_{1}.\ldots,X_{n})$. Then the distribution function of $X_{(n)}$ is 
\begin{align}
 F_{n}(t)&=\mathbb{P}(X_{(n)}<t)\nonumber\\
 &=\mathbb{P}(X_{1}<t)\ldots\mathbb{P}(X_{n}<t)\nonumber\\
 &=(1-e^{-t})^{n},~~t>0,\label{exp-dist}
\end{align}
and its density is 
\begin{equation}\label{density2}
 f_{n}(t)=n(1-e^{-t})^{n-1}e^{-t},~~t>0.
\end{equation}
For $s>0$ and a positive integer $m$, let $T_{m}$ be a gamma $G(s,m)$ random variable with density 
\begin{equation*}
 f(x|s,m)=\frac{s^{m}}{\Gamma(m)}e^{-sx}x^{m-1},~~x>0.
\end{equation*}
Note $T_{1}=G(s,1)=Exp(s)$, an exponential random variable with parameter ${s}> 0$.\vspace{0.1cm}\\
Peterson \cite{peter} computed the probability $\mathbb{P}(T_{1}>X_{(n)})$ in two different ways as follows:\vspace{0.1cm}\\ 
{\bf(A)} Conditioning on $X_{(n)}$:
\begin{align*}
 \mathbb{P}(T_{1}>X_{(n)})&=\mathbb{E}(\mathbb{P}(T_{1}>X_{(n)}|X_{(n)}))\nonumber\\
 &=\mathbb{E}(e^{-sX_{(n)}}).
\end{align*}
To compute the above expectation, he stated and proved the result (\cite[Lemma 1]{peter})
\begin{equation}\label{lemma1}
 X_{(n)}\stackrel{d}{=}\sum_{j=1}^{n}Y_{j},
\end{equation}
where $Y_{1},\ldots, Y_{n}$ are independent and $Y_{j}\sim Exp(j),~1\leq j\leq n.$ Note here $X\stackrel{d}{=}Y$ means $X$ and $Y$ have identical distributions. The proof of \eqref{lemma1} given in \cite{peter} is non-trivial and uses (a) the memoryless property of the exponential distribution and (b) the fact that $X_{(1)}=\min(X_{1},\ldots,X_{n})$ is an $Exp(n)$ random variable.\vspace{0.1cm}\\
{\bf(B)} Conditioning on $T_1$:
\begin{equation}\label{tcond}
 \mathbb{P}(T_1>X_{(n)})=\mathbb{E}(\mathbb{P}(X_{(n)}<T_{1}|T_{1})),
\end{equation}
which can easily be evaluated as $T_{1}\sim Exp(1)$. Equating \eqref{lemma1} and \eqref{tcond} he obtained the binomial identity in \eqref{identity1}.

\section{A Statistical Approach}
Our approach is based on computing the Laplace transformation of $X_{(n)}$ in two different ways and equate them to get the result in \eqref{identity1}.\vspace{0.1cm}\\
{\bf(i)} Using distribution function $F_{n}(t):$\\
Note that 
\begin{align}
 \mathbb{E}(e^{-sX_{(n)}})&=\int_{0}^{\infty}e^{-st}dF_{n}(t)\nonumber\\
 &=s\int_{0}^{\infty}F_{n}(t)e^{-st}dt~~~~~\text{(integration by parts)}\nonumber\\
 &=s\int_{0}^{\infty}(1-e^{-t})^{n}e^{-st}dt~~~~~(\text{use } \eqref{exp-dist})\nonumber\\
 &=\sum_{k=0}^{n}(-1)^{k}\binom{n}{k}\int_{0}^{\infty}e^{-(k+s)t}dt\nonumber\\
 &=\sum_{k=0}^{n}(-1)^{k}\binom{n}{k}\left(\frac{s}{s+k}\right)\label{identity2}\\
 &=f(s),~~(\text{say}).\nonumber
\end{align}
{\bf (ii)} Using the density of $f_{n}(t)$:\\ Using \eqref{density2}, we get
\begin{equation*}
 \mathbb{E}(e^{-sX_{(n)}})=n\int_{0}^{\infty}e^{-(s+1)t}(1-e^{-t})^{n-1}dt.
\end{equation*}
Making the transformation $w=(1-e^{-t})$, we get
\begin{align*}
 \mathbb{E}(e^{-sX_{(n)}})&=n\int_{0}^{1}(1-w)^{s}w^{n-1}dw\nonumber\\
 &=nB(s+1,n),
\end{align*}
where $B(a,b)$ is the usual beta function.\\
Note 
\begin{align*}
 nB(s+1,n)&=n\frac{\Gamma(s+1)\Gamma(n)}{\Gamma(n+s+1)}\\
 &=\frac{n!}{(s+1)\ldots (s+n)}.
\end{align*}
Thus, we get 
\begin{align}
 \mathbb{E}(e^{-sX_{(n)}})&=\prod_{k=1}^{n}\left(\frac{k}{s+k}\right)\label{identity3}\\
 &=g(s)~~~\text{(say)}.\nonumber
\end{align}
Equating \eqref{identity2} and \eqref{identity3}, we get the binomial identity in \eqref{identity1}.
\section{Remarks and Discussions}
Clearly, our statistical approach is simpler, direct and does not use any other property of the exponential distribution. Several remarks are in order.\\

(i) Let $Y_{1},\ldots, Y_{n}$ be $n$ independent random variables, where $Y_{j}\sim Exp(j),~1\leq j\leq n.$ Then from \eqref{identity3},

\begin{align*}
 \mathbb{E}(e^{-sX_{(n)}})&=\prod_{k=1}^{n}\left(\frac{k}{s+k}\right)\nonumber \\
 &=\prod_{k=1}^{n}\mathbb{E}(e^{-sY_{k}})\nonumber \\
 &=\mathbb{E}(e^{-s\sum_{k=1}^{n}Y_{k}}),
\end{align*}
for all $s>0$. By the uniqueness of the Laplace transform, we get
\begin{equation*}
 X_{(n)}\stackrel{d}{=}\sum_{k=1}^{n}Y_{k}.
\end{equation*}
This proves rigorously Lemma 1 given in \cite{peter}.\\

(ii) Note the binomial identity in \eqref{identity1} follows from 
\begin{equation*}
 f(s)=g(s),~~~\forall~ s>0.
\end{equation*}
This implies, obviously,
\begin{equation}\label{diffidentity}
 f(s)-sf'(s)=g(s)-sg'(s),~~~\forall~ s>0,
\end{equation}
which leads to 
\begin{equation}\label{diffidentity2}
 \sum_{k=0}^{n}\binom{n}{k}(-1)^{k}\left(\frac{s}{s+k}\right)^{2}=\prod_{k=1}^{n}\left(\frac{k}{s+k}\right)\left(\sum_{j=0}^{n}\frac{s}{s+j}\right).
\end{equation}
A proof of this result is given in \cite{peter} by computing the probability $\mathbb{P}(T_{2}>X_{(n)})$, where $T_{2}$ is a $G(s,2)$ random variable.\\

(iii) A curious reader might wonder if there exists a connection between these two approaches. Indeed, there is an interesting relationship between them, as seen follows:\\
Consider the general case 
\begin{align}
 \mathbb{P}(T_{m}>X_{(n)})&=\mathbb{E}(\mathbb{P}(T_{m}>X_{(n)})|X_{(n)})\nonumber\\
 &=\mathbb{E}\left(\sum_{k=0}^{m-1}e^{-sX_{(n)}}\frac{s^{k}}{k!}X_{(n)}^{k}\right)\nonumber\\
 &=\sum_{k=0}^{m-1}\frac{s^{k}}{k!}\mathbb{E}\left(e^{-sX_{(n)}}X_{(n)}^{k}\right),\label{identity4}
\end{align}
where we have used the fact (see \cite[equation (3.3.9)]{casella})
\begin{equation*}
 \mathbb{P}(T_{m}>x)=\sum_{k=0}^{m-1}\frac{e^{-sx}(sx)^{k}}{k!}.
\end{equation*}
Note $f(s)=\mathbb{E}(e^{-sX_{(n)}})$ has k-$th$ derivative
\begin{equation*}
 f^{(k)}(s)=(-1)^{k}\mathbb{E}\left(e^{-sX_{(n)}}X_{(n)}^{k}\right)
\end{equation*}
which, when substituted in \eqref{identity4}, leads to 
\begin{align}
 \mathbb{P}(T_{m}>X_{(n)})&=\sum_{k=0}^{m-1}(-1)^{k}\frac{s^{k}}{k!}f^{(k)}(s)\label{newidentity}\\
 &=\sum_{k=0}^{m-1}(-1)^{k}\frac{s^{k}}{k!}g^{(k)}(s)~~(\because f(s)=g(s))\label{newidentity2}
\end{align}
where $f$ and $g$ are defined respectively in \eqref{identity2} and \eqref{identity3}.

\noindent Note the equalities of the right-hand sides of \eqref{newidentity}  and \eqref{newidentity2} lead to binomial identities, while  \eqref{newidentity} (or \eqref{newidentity2}) indicates the connection between associated probabilities and the Laplace transforms.
For example, the case $m=2$ leads to \eqref{diffidentity} which in turn yields the identity in \eqref{diffidentity2}.\\
Note however $f'(s)$ (or $g'(s)$) is not of the form of the right side of \eqref{newidentity} (or \eqref{newidentity2}) and so it may not correspond to probability of some event involving exponential random variables. Hence, the binomial identity
\begin{equation*}
\prod_{k=1}^{n}\left(\frac{k}{s+k}\right)\sum_{j=1}^{k}\left(\frac{1}{j+s}\right)=\sum_{k=0}^{n}(-1)^{k+1}\binom{n}{k}\frac{k}{(k+s)^{2}},
\end{equation*}
which corresponds to $f'(s)=g'(s)$, may not have probabilistic interpretations. \\
 
(iv) For a general $m$, we get from right-hand sides of \eqref{newidentity} and \eqref{newidentity2} (or by directly computing $\mathbb{P}(T_{m}>X_{(n)})$),
\begin{equation*}
 \sum_{k=0}^{n}(-1)^{k}\binom{n}{k}\left(\frac{s}{s+k}\right)^{m}=\frac{n}{s}\sum_{k=0}^{m-1}\sum_{j=0}^{n-1}(-1)^{j}\binom{n-1}{j}\left(\frac{s}{s+j+1}\right)^{k+1},
\end{equation*}
which is different from the last identity given in \cite{peter}.
 Thus, several binomial identities can be obtained from the basic identity  in \eqref{identity1} and using \eqref{newidentity}  and \eqref{newidentity2}, all of which have probabilistic connections.\\ 

(v) Also, applying the well-known binomial inversion to the identity in \eqref{identity1}, we obtain for $s>0$,

\begin{equation}
 \sum_{k=0}^{n}(-1)^{k}\binom{n}{k}\prod_{j=1}^{k}\left(\frac{j}{s+j}\right)=\left(\frac{s}{s+n}\right),~~n\geq 0,
\end{equation}
where $\prod_{j=1}^{0}a_{j}$ is defined as unity. Similarly, from \eqref{diffidentity2}, we get

\begin{equation}
 \sum_{k=0}^{n}(-1)^{k}\binom{n}{k}\prod_{j=1}^{k}\left(\frac{j}{s+j}\right)\left(\sum_{i=0}^{k}\frac{s}{s+i}\right)=\left(\frac{s}{s+n}\right)^{2},
\end{equation}
and so on.\\ 
 Finally,  we hope the interested readers will enjoy these connections among binomial identities, probabilities of certain events and the Laplace transforms (and its derivatives) of  random variables associated with exponential distributions.
\printbibliography
\end{document}